\documentclass[12pt]{amsbook}
\usepackage{amssymb}
\usepackage{graphicx}
\usepackage{bbm}
\usepackage{float}
\usepackage{caption}
\usepackage{subcaption}
\newcommand{\T}{\mathbb{T}}
\renewcommand{\S}{\mathbb{S}}
\newcommand{\R}{\mathbb{R}}
\newcommand{\Z}{\mathbb{Z}}
\newcommand{\Q}{\mathbb{Q}}
\newcommand{\N}{\mathbb{N}}
\newcommand{\cR}{\mathcal R}
\newcommand{\cB}{\mathcal B}
\newcommand{\cD}{\mathcal D}
\newcommand{\cS}{\mathcal S}
\newcommand{\cT}{\mathcal T}
\newcommand{\tto}{\Rightarrow}
\newcommand{\transv}{\pitchfork}
\newcommand{\ep}{\begin{flushright}\fbox{}\end{flushright}}
\renewcommand{\le}{\mathop{\rm{LE}}}
\newcommand{\supp}{\mathop{\rm{supp}}}
\newcommand{\nuh}{\mathop{\rm{Nuh}}}
\newcommand{\pb}{\mathop{\rm{PB}}}
\newcommand{\per}{\mathop{\rm{Per}}}
\newcommand{\ehc}{{\rm{Ehc}}}
\newcommand{\eps}{\varepsilon}
\newcommand{\zeroeq}{\stackrel{\rm{o}}{=}}
\newtheorem{theorem}{Theorem}[chapter]

\newtheorem{proposition}[theorem]{Proposition}
\newtheorem{corollary}[theorem]{Corollary}
\theoremstyle{remark}
\newtheorem{remark}[theorem]{Remark}

\newtheorem{example}[theorem]{Example}

\numberwithin{equation}{chapter}

\title{A crash course in Pesin Theory}
\author{Jana Rodriguez Hertz}
\address{Department of Mathematics, SUSTech - China}
\address{SUSTech International Mathematical Center}
\email{rhertz@sustech.edu.cn}
\begin{document}
 \maketitle

\begin{figure}
\begin{center}
 \includegraphics[width=0.75\textwidth]{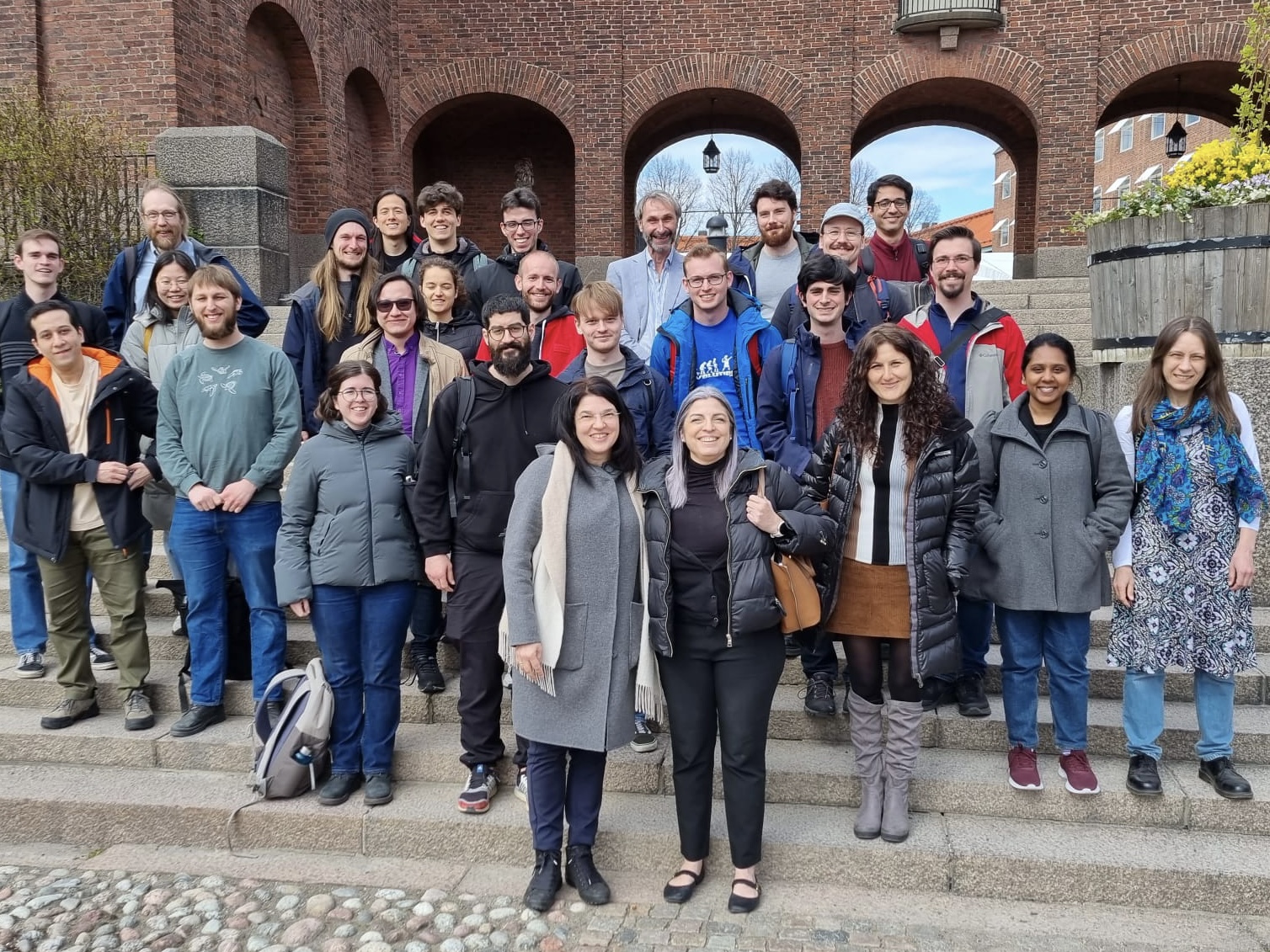}
 \caption*{Some of the people who attended the program ``Master class in low-dimensional dynamics'' at KTH, May 1-5, 2023.}
\end{center}
\end{figure}

\begin{figure}[h]
\begin{center}
 \includegraphics[width=0.75\textwidth]{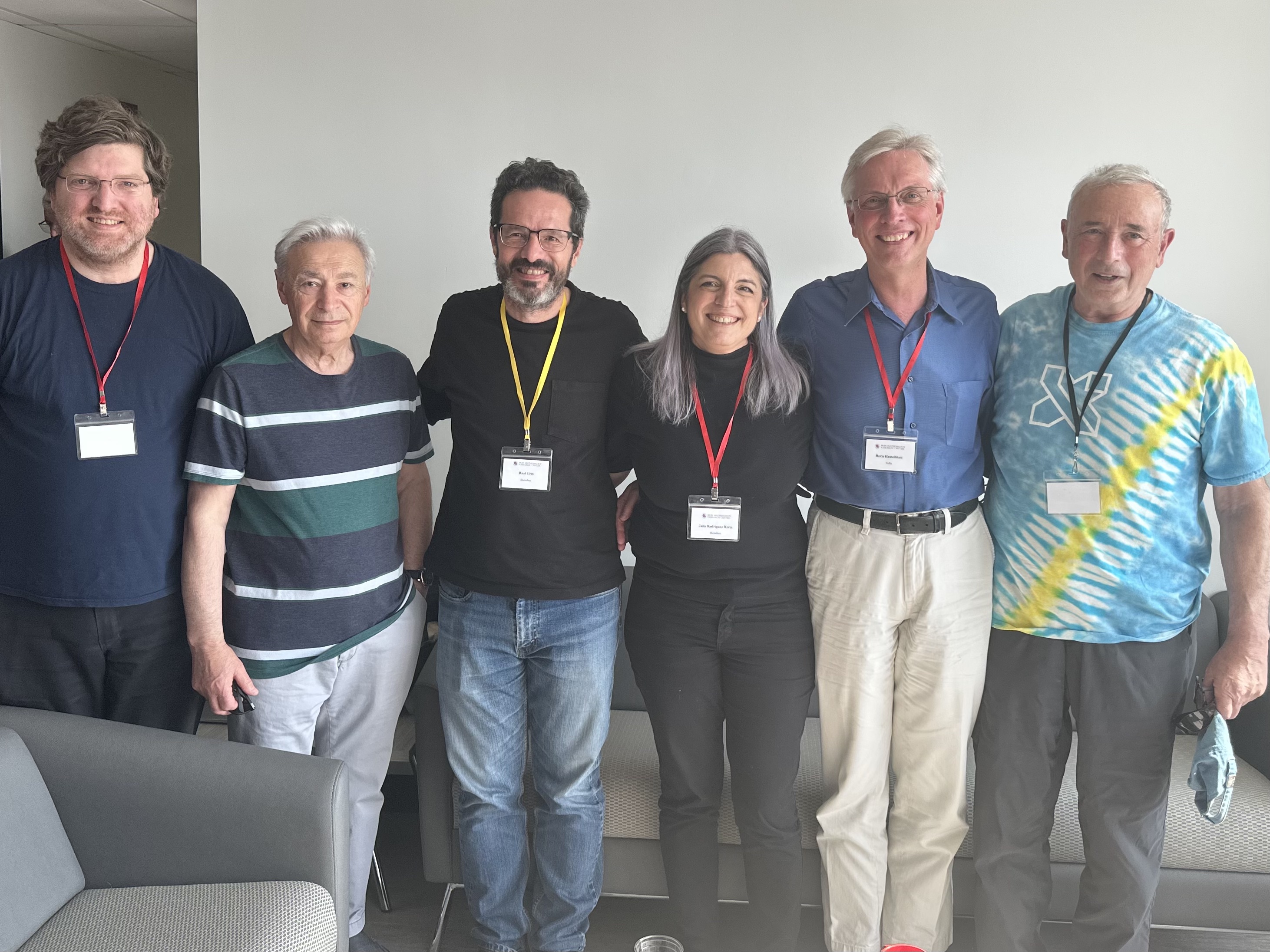}
 \caption*{Some of the people in the references, from left to right: Federico Rodriguez Hertz, Yasha Pesin, Ra\'ul Ures, me, Boris Hasselblatt, and Misha Brin}
\end{center}
\end{figure}

 \chapter*{  \dedicatory{\it To Yasha Pesin, Federico Rodriguez Hertz and Ra\'ul Ures, \\ in gratitude for their constant support\newline\newline\newline}
Disclaimer}

 This is more or less the content of a 5-hour mini-course that I have given at KTH in the program ``Master class in low-dimensional dynamics'' during May 1-5, 2023. It was not intended to be exhaustive but to give the students a grasp of basic ideas in Pesin Theory in a very short period. It is to be taken as a first approximation and as an invitation to study it in more detail. As such, there was a lot of hand-waving, probably there were mistakes and inaccuracies too. If you have any remark, correction, comment, you can contact me at {\tt janarhertz@gmail.com}. \par
 If you want to have a more complete approximation to Pesin Theory, I recommend the book by Luis Barreira and Yasha Pesin ``Introduction to Smooth Ergodic Theory''\cite{barreirapesin13} . Also a nice and friendly introduction to Pesin Theory is Mark Pollicott's ``Lectures on Ergodic Theory and Pesin Theory on compact manifolds'' \cite{pollicott}. You will have fun with any of these books. I have just learned from Sheldon Newhouse that Yuri Lima gave an excellent 5 day course on Pesin Theory that can be found on YouTube \cite{yurilima}. Yuri is a great expositor, so you should probably check on that too. \newline\par
 I would like to specially thank Kristian Bjerkl\"ov (KTH), Danijela Damjanovic (KTH), Liviana Palmisano (KTH), Maria Saprykina (KTH), Alan Sola (SU) and Michael Benedicks for warm hospitality during my stay. I would like to thank Rapha\"el Krikorian and Marco Martens for their Master Classes in the same school. Last but not least, I would like to thank the students for their great enthusiasm and very interesting questions and comments and making this school fun. \newline\par
 \thanks{JRH  was supported by NSFC 12161141002 and NSFC 12250710130}

 \chapter{Leaving hyperbolicity}\label{hyperbolicity}
 \section{Arnold's cat map}\label{Arnold}
 Let us begin by looking at (uniform) hyperbolic behavior. This will be our toy example:
 $$A=\left(
\begin{array}{cc}
2&1\\1&1 
\end{array}
\right)$$
We will consider this linear automorphism acting in the 2-torus $\T^{2}=\R^{2}/\Z^{2}$:
$$f(x)=Ax\qquad \mod \Z^{2}.$$
There exists $\mu>1$ such that $\mu$ and $\mu^{-1}$ are the eigenvalues of $A$ (why?). Let $v_{+}$ and $v_{-}$ be the corresponding eigenvectors. Then, under successive iterations $Df(x)$ stretches and contracts vectors exponentially, respectively. 
$$Df^{n}(x)v_{\pm}=\mu^{\pm n}v_{\pm}\qquad \forall n\in\N_{0}.$$
In the future, we will want to identify the rate of expansion of vectors. Assume that $||v_{\pm}||=1$. Then 
the average rate of expansion is 
\begin{equation}\label{rateexpansion1}
||Df^{n}(x)v_{\pm}||^{\frac{1}{n}}=\mu^{\pm 1}. 
\end{equation}

In this case, we have eigenvalues, so it is easy to compute the average rate of expansion. In general, we will not have an invariant direction, and even if we had, we will not have a fixed rate of expansion. We want to find out who plays the role of $\mu$ in general.
\subsection{Exercises - Section \ref{Arnold}}
\begin{enumerate}
 \item Check that $f:\T^{2}\to \T^{2}$ is well defined.
 \item Check that periodic points of $f$ are dense in $\T^{2}$. Prove that $p\in\T^{2}$ is periodic if and only if $p\in\Q^{2}/\Z^{2}$. 
 \item Compute $\mu,\mu^{-1}$ the eigenvalues of $A$ and $v_{+},v_{-}$, the eigenvectors. 
\end{enumerate}
\chapter{Entering non-uniform hyperbolicity}\label{NUH}
\section{Lyapunov exponents}\label{section.LE}
Remember we want to estimate the average expansion rate. Since in general $\mu$ in formula \eqref{rateexpansion1} is not independent of $n$, we could just settle for the asymptotic average expansion rate. That is, we could find out, for any given $v\in T_{x}M$, who is:
$$\limsup_{n\to\infty}||Df^{n}(x)v||^{\frac{1}{n}}.$$
For technical and historical reasons, which will be clear later (or not), we choose instead to compute the logarithm of this rate. That is for any $x\in M$ and $v\in T_{x}M\setminus\{0\}$, 
\begin{equation}\label{lyapunovexponents}
 \le(x,v)=\limsup_{n\to\infty}\frac{1}{n}\log||Df^{n}(x)v||.
\end{equation}
We will call this amount the {\em Lyapunov exponent of $x$ in the direction of $v$}.\par
In the case of Arnold's cat map we have:
$$\le(x,v_{+})=\log \mu>0\quad\text{and}\qquad \le(x,v_{-})=-\log\mu<0$$
A positive Lyapunov exponent means average exponential expansion in the direction of $v$. A negative Lyapunov exponent means average exponential contraction in the direction of $v$. Zero Lyapunov exponent means that whatever contraction or expansion there is, it is subexponential. It is not seen with the exponential glasses. \par
Note that the Lyapunov exponents does not take into account the modulus of $v$, it only depends on the direction of $v$. That is, if $v\in T_{x}M\setminus\{0\}$ and $\alpha\ne0$, 
$$\le(x,\alpha v)=\le(x,v). $$
Can you prove this?\newline \par
What about $\le(x,u+v)$ for $u,v\in T_{x}M\setminus\{0\}$? Let $$K=\max(\le(x,u),\le(x,v)).$$
Then, 
$$
 \frac{1}{n}\log||Df^{n}(x)(u+v)||\leq \frac{1}{n}\log (||Df^{n}(x)u||+||Df^{n}(x)v||)$$
 Since $\max(\le(x,u),\le(x,v))\leq K$, for large $n\in\N$ and small $\eps>0$, 
 $$\max\left(\frac{1}{n}\log ||Df^{n}(x)u||,\frac{1}{n}\log ||Df^{n}(x)u|\right)|\leq  K+\eps, $$
 so we have 
 $$ \frac{1}{n}\log||Df^{n}(x)(u+v)||\leq \frac{1}{n}\log 2\exp(n(K+\eps))=\frac{\log 2}{n}+K+\eps.$$
 Since this inequality holds for each $\eps>0$ for all sufficiently large $n\in \N$, we have
 $$\le(x,u+v)\leq\max(\le(x,u),\le(x,v)).$$
If $\le(x,u)<\le(x,v)$, then
\begin{eqnarray*}
\le(x,v)&=&\le(x,u+v-u)\leq\max(\le(x,u+v),\le(x,u))\\
&\leq&\max(\le(x,u+v),\le(x,u))=\le(x,u+v)\leq\le(x,v) 
\end{eqnarray*}
What happens if $\le(x,u)=\le(x,v)$?
\begin{proposition}\label{prop.LE} Let $f:M\to M$ be a diffeomorphism. Let $u,v\in T_{x} M\setminus{0}$ with $x\in M$.
\begin{enumerate}
 \item $\le(x,\alpha v)=\le(x,v)$ for all $\alpha\ne0$. 
 \item $\le(f(x),Df(x)v)=\le(x,v)$.
\item $\le(x,u+v)\leq\max(\le(x,u),\le(x,v)).$
\item if $\le(x,u)<\le(x,v)$ then $\le(x,u+v)=\le(x,v)$. 
\end{enumerate}
\end{proposition}
\subsection{Exercises - Section \ref{NUH}.\ref{section.LE}}
\begin{enumerate}
 \item Can you compute {\em all} Lyapunov exponents for the Arnold's cat map? (Section \ref{hyperbolicity}.\ref{Arnold})
 \item Check items 1 and 2 of Proposition \ref{prop.LE}.
 \item Prove that the inequality in item 3 of Proposition \ref{prop.LE} can be strict.  
 \item Check that the Lyapunov exponents are constant $\mu$-almost everywhere if $\mu$ is an ergodic invariant measure. 
\end{enumerate}
\section{Oseledets theorem}\label{section.oseledets}
Let us review the notion of hyperbolicity: given a $C^{1}$ diffeomorphism $f:M\to M$ a compact $f$-invariant set $\Lambda$ is {\em hyperbolic} if there is a $Df$-invariant splitting of $TM$ over $\Lambda$, $T_{x}M=E_{x}^{s}\oplus E^{u}_{x}$, with $x\in M$, that is $Df(x)E^{*}_{x}=E^{*}_{f(x)}$ with $*=s,u$ for all $x\in \Lambda$, and constants $C>0$ and $0<\lambda<1$ such that for all unit vectors $v^{s}\in E^{s}$ and $v_{u}\in E^{u}$ it holds:
\begin{enumerate}
 \item $\|Df^{n}(x)v^{s}\|\leq C\lambda^{n}$ for all $n\geq 0$,
 \item $\|Df^{-n}(x)v^{u}\|\leq C\lambda^{n}$ for all $n\geq 0$.
\end{enumerate}
Equivalently, if there exists a Riemannian metric (called {\em adapted metric}) such that for all unit vectors
$v^{s}\in E^{s}$ and $v_{u}\in E^{u}$ it holds:
$$\|Df(x)v^{s}\|<1<\|Df(x)v^{u}\|.$$
(This is not a trivial fact. In a future version -with a coauthor- we will explain this better)\newline
What happens in the general case? Can there be infinitely many Lyapunov exponents for $x\in M$? The answer is no. You will prove it in the set of exercises. What is the role of Lyapunov exponents? Are there invariant directions? Can we expect something similar to hyperbolicity?
The following theorem partially answers the previous questions:
\begin{theorem}[Oseledets theorem]\label{theorem.oseledets}
Given a $C^{1}$ diffeomorphism $f:M\to M$, for a total measure set $\cR$ (that is, a set that satisfies $\mu(\cR)=1$ for all $f$-invariant measure $\mu$) for all $\eps>0$ there is a Borel function $C_{\eps}:\cR\to(1,\infty)$ such that for all $x\in \R$:
\begin{enumerate}
 \item $T_{x} M=E_{\lambda_{1}(x)}(x)\oplus \dots\oplus E_{\lambda_{k_{x}}(x)}(x)$ Oseledets splitting, such that $LE(x,v)=\lambda_{i}(x)$ for all $v\in E_{\lambda_{i}(x)}(x)\setminus\{0\}$.
 \item For all unit vector $v\in E_{\lambda}(x)$ and all $n\in \Z$:
 $$C_{\eps}^{-1}(x)\exp(\lambda-\eps)n\leq\|Df^{n}(x)v\|\leq C_{\eps}(x)\exp(\lambda+\eps)n,$$
 \item $\angle(E_{\lambda}(x),E_{\lambda'}(x))\geq C_{\eps}^{-1}(x)$ if $\lambda\ne \lambda'$,
 \item $$\exp(-\eps)\leq \frac{C_{\eps}(f(x))}{C_{\eps}(x)}\leq \exp \eps$$
\end{enumerate} 
\end{theorem}
If you look carefully at the Oseledets theorem above, it sort of behaves like hyperbolicity when the Lyapunov exponent is different from zero, except that the constant $C>0$ in the definition of hyperbolicity now depends on $x$. However, its decay with the action of $f$ is tempered. \par
Note that the set $\cR$ could be in general pretty small from the topological point of view. If we consider, for instance, the North pole-South pole dynamics $f$ on the circle $\S^{1}$, we get a $2$ point $\cR$.  Indeed, for $f$ we have two hyperbolic fixed points $N$ and $S$ such that all orbits except $N$ converge to $S$ in the future, and all orbits except $S$ converge to $N$ in the past. Then $\cR=\{N,S\}$.\par
However, in this course we will be mainly considering volume preserving diffeomorphisms $f$, so $\cR$ will be ``large'' in our setting. \newline\par
Many times we will be mainly concerned about the sign of the Lyapunov exponents. To that end, we define the {\em zipped Oseledets splitting} like this:
$$T_{x}M= E^{-}_{x}\oplus E^{0}_{x}\oplus E^{+}_{x}\quad \forall x\in \cR$$
where
\begin{itemize}
 \item $\le(x,v)<0$ for all $v\in E^{-}_{x}$,
 \item $\le(x,v)=0$ for all $v\in E^{0}_{x}$,
  \item $\le(x,v)>0$ for all $v\in E^{+}_{x}$
\end{itemize}
If $\mu$ is an invariant measure such that $E^{0}_{x}=\{0\}$ for $\mu$-almost every $x\in M$, then $\mu$ is a {\em hyperbolic measure}. If $f$ is a volume preserving diffeomorphism and the volume measure $m$ is hyperbolic, we say $f$ is {\em non-uniformly hyperbolic}. \newline\par
If we fix $\eps>0$ and $L>1$, and consider the points $x\in M$ such that $C_{\eps}(x)\leq L$ we obtain ``hyperbolicity'', but we lose invariance. These sets are important in Pesin theory, and are called {\em Pesin blocs}. Namely,
$$
\pb(\eps,L)=\{x\in M: C_{\eps}(x)\leq L\}. 
$$

For any fixed $\eps>0$ we have
\begin{equation}\label{pesin.block}
\cR=\bigcup_{L\in\N}\pb(\eps,L). 
\end{equation}

We will be back on Pesin blocs later. But for the time being, it is important to remember that in the non-uniformly hyperbolic setting there is a trade-off between the uniformity of hyperbolicity and invariance. \newline
\begin{center}
\begin{tabular}{ccc}
hyperbolicity &$\leftrightarrow$ & invariance\\
+ uniform &$\to$&-invariant\\
-uniform&$\leftarrow$&+invariant
\end{tabular}
\end{center}
\vspace*{1em}
There is a trick that sometimes is helpful to achieve both: let us remember Poincar\'e Recurrence Theorem:
\begin{theorem}[Poincar\'e Recurrence Theorem] If $f: M\to M$ preserves a probability measure $\mu$, then for all measurable sets $B$ such that $\mu(B)>0$, for $\mu$-almost every $x\in B$ there exist infinitely many $n\in \N$ such that $f^{n}(x)\in B$. 
\end{theorem}
Because of Equation \eqref{pesin.block}, for all invariant probability measure, there is $L>0$ such that $\mu(\pb(\eps, L))>0$. So we have sort of hyperbolic iterations. 

\subsection{Exercises}
\begin{enumerate}
\item If $\Lambda$ is a hyperbolic set, check that the Lyapunov exponents of $x\in \Lambda$ corresponding to vectors in $E^{s}$ are negative, and corresponding to vectors in $E^{u}$ are positive. Check that there are no vectors with zero Lyapunov exponents.
 \item Prove that there cannot be more than $n$ different Lyapunov exponents for $x$, where $n=\dim M$. 
 \item Prove that $f(\pb(\eps,L))\subset \pb(\eps,(\exp \eps)L)$
\end{enumerate}
\section{Other weak forms of hyperbolicity}\label{section.weak.hyperbolicity}
For future use, let us state here other weak forms of hyperbolicity. One is partial hyperbolicity. We say a diffeomorphism $f:M\to M$ is {\em partially hyperbolic} if there is a $Df$-invariant splitting of the tangent bundle $T_{x}M=E^{s}_{x}\oplus E^{c}_{x}\oplus E^{u}_{x}$ and a Riemannian metric such that for all unit vectors $v^{*}\in E^{*}_{x}$ with $*=s,c,u$ we have:
$$\|Df(x)v^{s}\|<1<\|Df(x)v^{u}\|\qquad\text{and}$$
$$\|Df(x)v^{s}\|<\|Df(x)v^{c}\|<\|Df(x)v^{u}\|.$$
An example of a partially hyperbolic non hyperbolic diffeomorphism is $f:\T^{3}\to\T^{3}$ such that 
$f=A\times {\rm id}$, where $A$ is as in Section \ref{Arnold} and ${\rm id}$ is the identity on the circle $\T^{1}$. \newline\par
Yet another example of a weaker form of hyperbolicity is the following: when there is not necessarily expansion or contraction, but there is an invariant splitting such that $Df$ contracts stronger in one of the bundles than in the other. We call this $f$ has a {\em dominated splitting} Namely, there is a $Df$-invariant splitting $TM=E\oplus F$ such that 
for each unit vector $v_{E}\in E$ and $v_{F}\in F$ we have:
$$\|Df(x)v_{E}\|<\|Df(x)v_{F}\|.$$
Clearly, if $f$ is a partially hyperbolic diffeomorphism, then it has a dominated splitting. 
\subsection{Exercises}
\begin{enumerate}
 \item Prove that if there is a a dominated splitting $TM=E\oplus F$, then $\le(x,v_{E})<\le(x,v_{F})$ for all unit vectors $v_{E}\in E$ and $v_{F}\in F$.
 \item Also, if there is a dominated splitting, for all $x\in \cR$ in the Oseledets Theorem, there is an $i(x)$ such that $E_{\lambda_{1}}(x)\oplus\dots\oplus E_{\lambda_{i(x)}}(x)=E_{x}$ and 
 $E_{\lambda_{i(x)+1}}\oplus\dots\oplus E_{\lambda_{k(x)}}= F_{x}$. 
\end{enumerate}
\section{Our first true example of non-uniform hyperbolicity}
Let us see a non-hyperbolic non-uniformly hyperbolic example dynamical system. Let us go back to our favorite toy example introduced in Section \ref{Arnold}. 
$$A=\left(
\begin{array}{cc}
2&1\\
1&1 
\end{array}
\right)$$
Then, for $x\in \T^{2}$ we have:
$$A(-x)=-Ax.$$
If we quotient $\T^{2}$ by the equivalent relation $x\sim -x$, we obtain the sphere $\S^{2}$! You can check this by using a piece of paper (improved picture later): 
\begin{figure}[h]
\begin{center}
 \includegraphics[width=\textwidth]{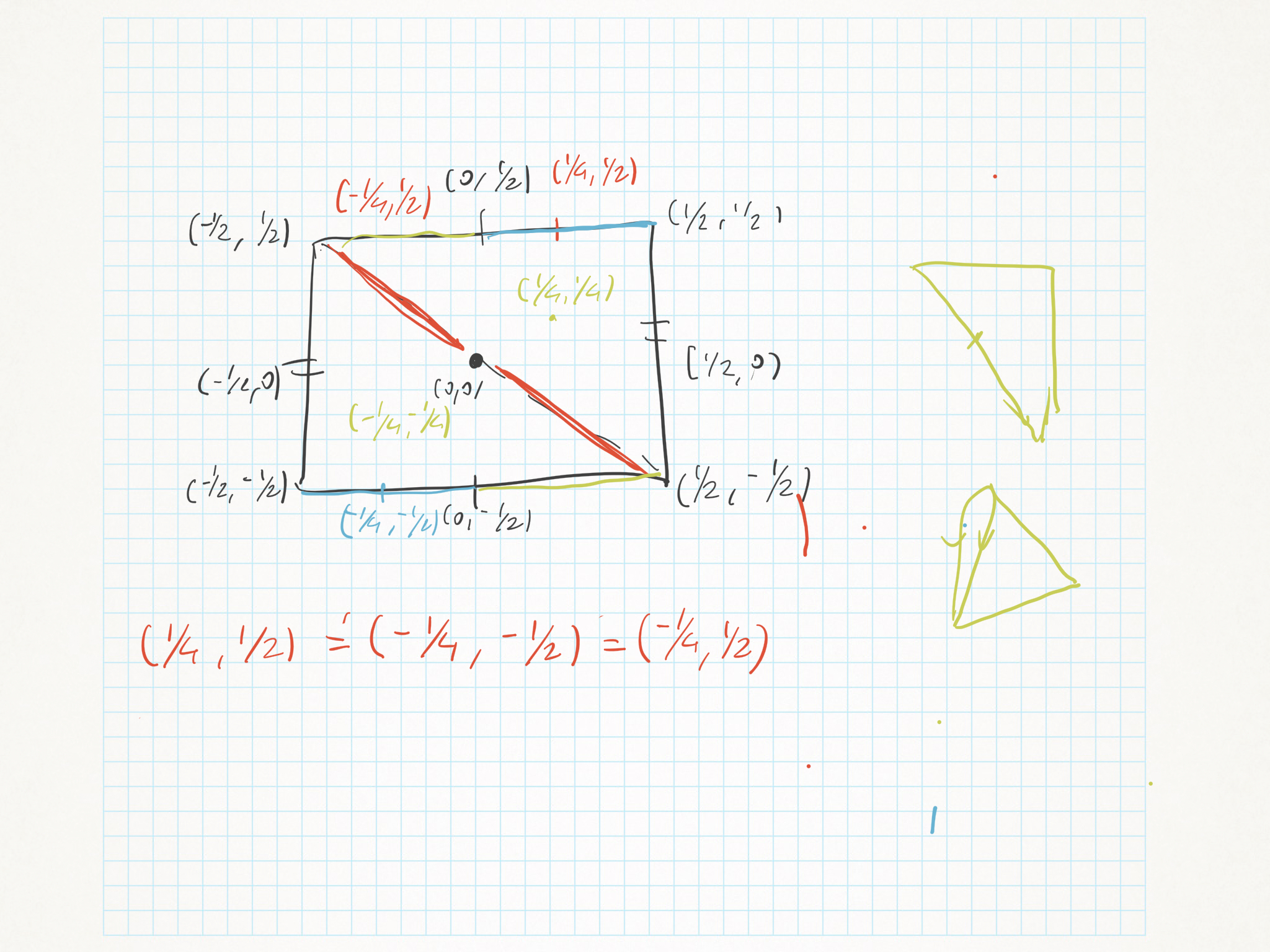}
\end{center}
\end{figure}
There is a projection $\pi:\T^{2}\to \S^{2}$, such that $\pi^{-1}(x)=\{x,-x\}$, and this projection makes the following diagram commute:
$$
\begin{array}{rcccl}
 &&A&&\\
& \T^{2}&\to&\T^{2}&\\
\pi& \downarrow&&\downarrow&\pi\\
& \S^{2}&\to&\S^{2}&\\
& &\tilde f&&
\end{array}
$$
where $\tilde f$ is the induced map by $A$ and $\pi$, namely: $\tilde f(\pi(x))=\pi(Ax)$. \par
Let $m_{\T^{2}}$ be the Lebesgue measure in the 2-torus. Then, the normalized Lebesgue measure on the 2-sphere is, for any measurable set $B$:
$$m_{\S^{2}}(B)=m_{\T^{2}}(\pi^{-1}(B))$$
We call this example the {\em pseudo-Anosov map of the sphere}. There are other examples of pseudo-Anosov maps in arbitrary compact surfaces. These are also non-uniformly hyperbolic. There are also $C^{0}$-perturbation maps by Katok, that are non-uniformly hyperbolic, now known as {\em Katok maps} \cite{katok79}. 
\subsection{Exercises}
\begin{enumerate}
 \item Prove that $\tilde f$ preserves $m_{\S^{2}}$.
 \item Can you prove that $\tilde f$ is non-uniformly hyperbolic?
\end{enumerate}
\chapter{Absolute continuity}
\section{Stable and unstable manifolds}
We are interested in identifying points that go asymptotically to the same places as $x$ in the future and in the past. In (uniformly) hyperbolic dynamics, we usually define the {\em stable manifold} by
$$W^{s}(x)=\{y\in M: d(f^{n}(x),f^{n}(y))\to 0 \text{ as }n\to\infty\},$$
and the {\em unstable manifold} by 
$$W^{u}(x)=\{y\in M: d(f^{-n}(x),f^{-n}(y))\to 0 \text{ as }n\to\infty\}.$$
In non-uniformly hyperbolic dynamics, we will want to be more specific. We are interested in points that asymptotically get closer at an exponential rate. We call these the {\em Pesin stable manifolds} 
 $$W^{-}(x)=\left\{y\in M:\limsup_{n\to \infty}\frac{1}{n}\log d(f^{n}(x),f^{n}(y))<0\right\}$$
 and the {\em Pesin unstable manifolds} will be
  $$W^{+}(x)=\left\{y\in M:\limsup_{n\to \infty}\frac{1}{n}\log d(f^{-n}(x),f^{-n}(y))<0\right\}$$
  Sometimes we will be interested in the local versions of these sets, namely, the {\em local Pesin stable manifolds}
  $$W^{-}_{\eps}(x)=\left\{y\in W^{-}(x): d(f^{n}(x),f^{n}(y))\leq \eps\quad \forall n\geq 0\right\},$$
  and the {\em local Pesin stable manifold} by
    $$W^{+}_{\eps}(x)=\left\{y\in W^{+}(x): d(f^{-n}(x),f^{-n}(y))\leq \eps\quad \forall n\geq 0\right\}.$$
    When we want to emphasize that we are considering local Pesin manifolds, but the size $\eps>0$ is not relevant, we will denote $W^{-}_{loc}(x)$ and $W^{+}_{loc}(x)$. \newline\par
    We can recall from hyperbolic dynamics the following theorem   
\begin{theorem}[Stable manifold theorem \cite{HPS77}]
Let $f:M\to M$ be a $C^{1}$-diffeomorphism, and $\Lambda$ a compact invariant set that is hyperbolic for $f$.  Then for all $x\in\Lambda$, $W^{s}_{loc}(x)$ contains an $s$-disc centered at $x$, where $s=\dim E^{s}_{x}$, such that 
\begin{enumerate}
 \item $T_{y}W^{s}_{loc}(x)=E^{s}_{y}$ for all $y\in W^{s}_{loc}(x)$, and
 \item $x\mapsto W^{s}_{loc}(x)$ continuous in the $C^{1}$-topology for all $x\in \Lambda$.
\end{enumerate}
\end{theorem}
Of course, we have an analogous statement for local unstable manifolds. More surprisingly, we have an analogous statement for local Pesin stable manifolds. 
\begin{theorem}[Pesin stable manifold theorem \cite{pesin76}] Let $f:M\to M$ be a $C^{2}$-diffeomorphism. Let $\mu$ be any invariant probability measure. Then for all $\mu$- almost every $x\in\cR$, $W^{-}_{loc}(x)$ contains an $s$-disc centered at $x$, where $s=\dim E^{-}_{x}$, such that 
\begin{enumerate}
 \item $T_{y} W^{-}_{loc}(x)=E^{-}_{y}$ for all $y\in W^{-}_{loc}(x)$, 
 \item if $\mu$ is ergodic, then $x\mapsto W^{-}_{loc}(x)$ varies continuously on $\pb(\eps,L)$ wherever it is well-defined. 
\end{enumerate}
\end{theorem}
If $\mu$ is not ergodic, then we can also get item 2. above, but we have to consider the points in $\pb(\eps,L)$ such that the dimension of $E^{-}_{x}$ is constant.
\begin{remark} The $C^{2}$ hypothesis is essential in Pesin stable manifold theorem. Actually, asking $C^{1+\alpha}$ is enough. But $C^{1}$ is not. There is an example of a $C^{1}$ diffeomorphism by Pugh \cite{pugh84}, for which there is an invariant probability measure $\mu$ and a measurable set $B$ with $\mu(B)>0$ such that $W^{s}_{loc}(x)$ is not a disc for all $x\in B$.
\end{remark}
\begin{remark}
 However, the Pesin Stable Manifold Theorem holds for $C^{1}$-diffeomorphisms if in the zipped Oseledets splitting we have this:
 $$T_{x}M=E^{-}_{x}\oplus_{<}E^{0}_{x}\oplus E^{+},$$
that is, the bundle $E^{0}_{x}\oplus E^{+}_{x}$ dominates $E^{s}_{x}$ for all $x\in M$. This is due to Abdenur, Bonatti and Crovisier.  
\end{remark}
  \subsection{Exercises}
\begin{enumerate}
 \item Prove that for our favorite toy example $A=\left(
\begin{array}{cc}
 2&1\\1&1
\end{array}\right)
$, we have that $y\in W^{-}(x)$ if and only if $x-y\in \R v_{-}$, and $y\in W^{+}(x)$ if and only if $x-y\in \R v_{+}$. You can use either stable or Pesin stable manifold definition. 
\item Prove that for all $\eps>0$ small
$$W^{-}(x)=\bigcup_{n\geq0}f^{-n}(W^{-}_{\eps}(f^{n}(x)))\quad\text{and}\quad W^{+}(x)=\bigcup_{n\geq 0}f^{n}(W^{+}_{\eps}(f^{-n}(x))).$$
\item Prove that if $f$ is hyperbolic on $\Lambda$, then $W^{s}(x)=W^{-}(x)$ and $W^{u}(x)=W^{-}(x)$ for all $x\in \Lambda$. 
\end{enumerate}
\section{Absolute continuity}
Given a manifold $M$, we say that a partition $\xi$ is {\em measurable} if $M/\xi$ can be separated by a countable number of measurable sets. 
\begin{example}
 The partition $\xi$ of $\T^{2}$ by vertical circles is measurable. $(M/\xi)\sim \S^{1}$, which can be separated by the rationals. 
\end{example}
\begin{example}
 The partition of $\T^{2}$ by stable manifolds of $A$ (our example of Section \ref{Arnold}) is {\em not} measurable. Don't try to prove it, just feel it. 
\end{example}
Let us recall that a Lebesgue space is a space which is isomorphic to the closed interval $[0,1]$ with the Lebesgue measure union a finite (possibly empty) family of atoms. 
\begin{proposition}
 If $M$ is a Lebesgue space and $\xi$ is a measurable partition, then $M/\xi$ is a Lebesgue space. 
\end{proposition}
Let $(M,\mathcal{B}, m)$ be a Lebesgue space and $\xi$ be a measurable partition. Associated to these, there is a canonical system of {\em conditional measures}
$$x\mapsto m^{\xi}_{x}$$
such that $m^{\xi}_{x}$ is the (normalized) Lebesgue measure on $\xi(x)$, the element of $\xi$ that contains $x$. This system has the property that for all $B\in \mathcal{B}$ the set $B\cap \xi(x)$ is measurable in $\xi(x)$ for $m$-almost all $x\in M$, and the function $$x\mapsto m^{\xi}_{x}(B\cap\xi(x))$$
is measurable. Moreover, 
\begin{equation}\label{fubini}
m(B)=\int_{M}m^{\xi}_{x}(B\cap \xi (x))dm
\end{equation}
\begin{figure}[h]
\begin{center}
 \includegraphics[width=0.5\textwidth]{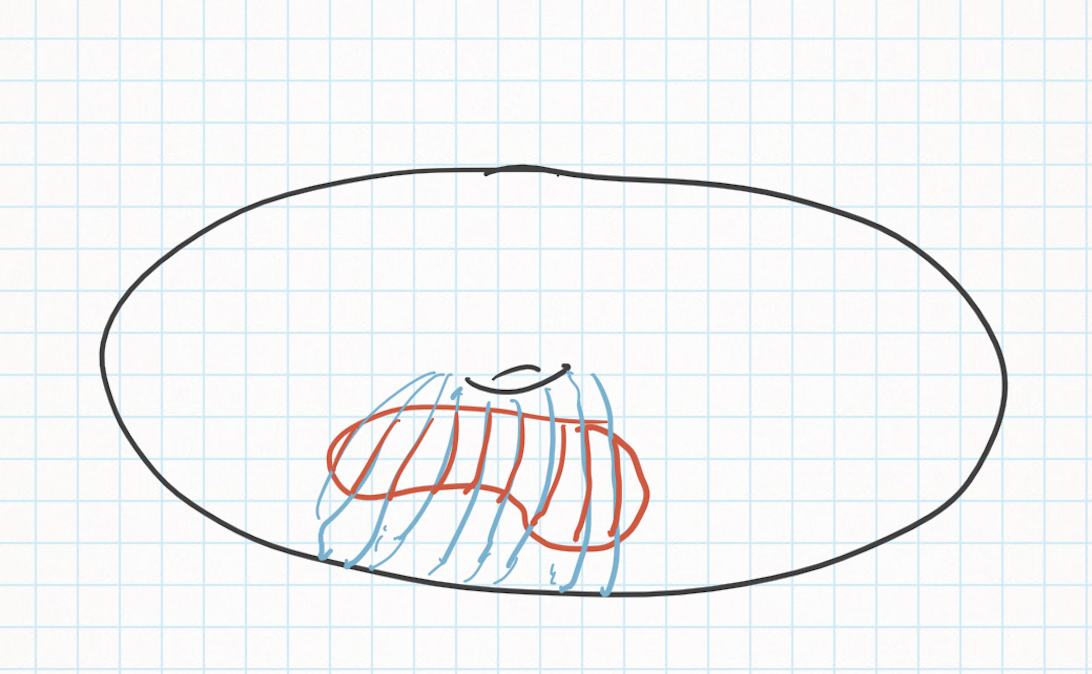}
 \caption{A set of conditional measures}
\end{center}
\end{figure}
For each measurable partition this canonical system of conditional measures is unique modulo a zero measure set. Any other system coincides with it for $m$-almost all $x\in M$. Conversely, if there is a canonical system for a partition, then the partition is measurable. \par
We are interested in stable and unstable partitions, but these partitions are not measurable in general. So, another method has to be devised.  \par
A measurable partition $\xi$ is {\em subordinate} to the (Pesin) unstable partition $W^{-}$ if for $m$-almost every $x$ we have $\xi(x)\subset W^{-}(x)$, and $\xi(x)$ contains a neighborhood of $x$ which is open in the topology of $W^{-}(x)$.  
\begin{theorem}[Pesin \cite{pesin76}] Let $f:M\to M$ be a $C^{2}$-diffeomorphism that preserves Lebesgue measure. Then there are partitions $\xi^{\pm}$ subordinate to $W^{\pm}$ respectively.
\end{theorem}
We say that the volume measure $m$ has {\em absolutely continuous conditional measures on unstable manifolds} if for all measurable partition $\xi$ subordinate to $W^{+}$, we have that
$$m^{\xi}_{x}<<Leb^{+}_{x}\quad \text{and}\quad  Leb^{+}_{x}<<m^{\xi}_{x}\quad \text{for }m-{a.e. }x, $$
where $Leb^{+}_{x}$ is the Riemannian measure on $W^{+}(x)$ given by the Riemannian structure of $W^{+}(x)$ inherited from $M$. Here $<<$ means absolutely continuous. Remember $\mu$ is absolutely continuous with respect to $\nu$ ($\mu<<\nu$) if for all measurable set $A$ we have $\nu(A)=0\tto \mu(A)=0$.\newline\par
Take $x_{0}\in M$ a Lebesgue density point of the set $B^{+}(u)=\{x\in\cR:\dim E^{+}_{x}=u\}$, that is $x$ is a point such that 
$$\lim_{\eps\to0}\frac{m(B_{\eps}(x_{0})\cap B^{+}(u))}{m(B_{\eps}(x_{0}))}=1.$$
(Remember that if $m(B^{+}(u))>0$ then $m$-almost every point in $x\in B^{+}(u)$ is a Lebesgue density point of $B^{+}(u)$. This is true for any measurable set). Now take two small discs $T,T'$ near $x_{0}$ transverse to $W^{+}_{loc}(x_{0})$, as in the figure:
\begin{figure}[h]
\begin{center}
 \includegraphics[width=\textwidth]{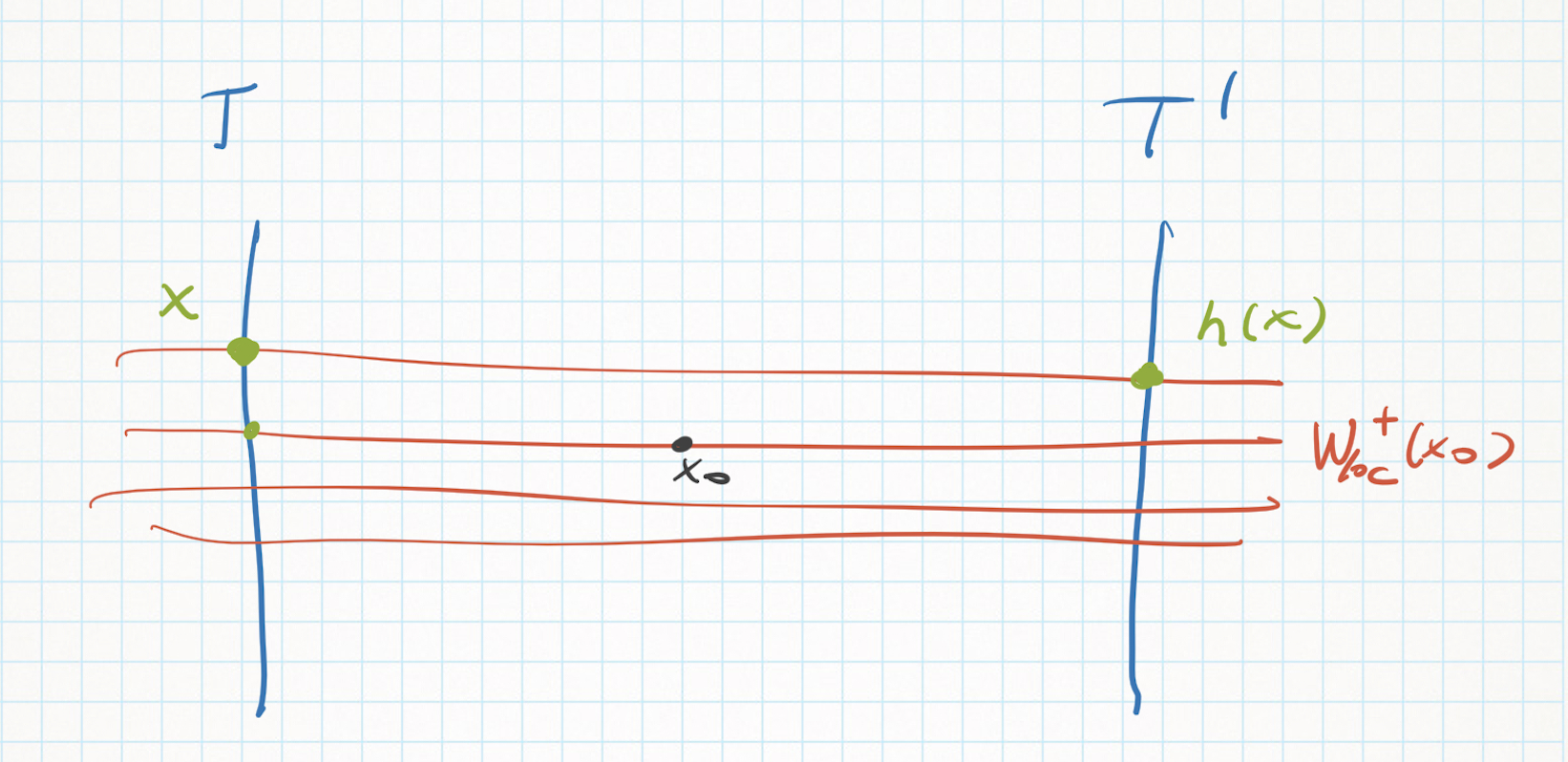}
 \caption{Unstable holonomy map}
\end{center}
\end{figure}
We will define the holonomy map with respect to these transversals by $h:D\subset T\to T'$
such that $$h(x)=W^{+}_{loc}(x)\cap T'.$$
The domain of $h$ is $D$ which is the set of points in $T\cap B(u)$ such that $W^{+}_{loc}(x)\cap T'\ne\emptyset$. $h$ is a bijection in its domain and image. \par
We say that the Pesin unstable partition is (transversely) {\em absolutely continuous} if all holonomy maps are measurable and take Lebesgue zero measure sets of $T$ into Lebesgue zero measure sets of $T'$. Analogous definition holds for the Pesin stable partition. 
\begin{theorem}[Pesin] \label{pesin.absolute.continuity}Let $f:M\to M$ be a $C^{2}$-diffeomorphism preserving a smooth volume. Then the Pesin stable and unstable partitions are absolutely continuous. 
\end{theorem}
\begin{remark} Let $$B^{+}(u)=\{x:\dim E^{+}_{x}=u\}.$$
Then, when restricted to the set $\pb(\eps,L)\cap B^{+}(u)$, the unstable holonomy maps are continuous. Analogous statement holds for the stable holonomy maps. 
\end{remark} 
\begin{remark} Pesin Theorem \ref{pesin.absolute.continuity} implies there exist measurable partitions subordinate to $W^{\pm}$ which have the Fubini-type property described above. In particular, they have absolutely continuous conditional measures. 
\end{remark}
\chapter{The Hopf argument}
\section{Ergodicity}
One of the main uses of Pesin Theory is getting tools to establish ergodicity. Remember that for a diffeomorphism $f:M\to M$ we say that an $f$-invariant probability measure $\mu$ is {\em ergodic} if any measurable set $A\in\cB$ that is invariant modulo a zero set, that is, $\mu(A\triangle f^{-1}(A))=0$ then we have either $\mu(A)=0$ or $\mu(A)=1$.\par
We will mention the following classical result, adapted to our setting
\begin{theorem}[Birkhoff Ergodic Theorem] Let $f:M\to M$ be a diffeomorphism preserving the probability measure $m$. Then for all continuous maps $\varphi:M\to \R$ (that is $\varphi\in C^{0}(M)$), we have:
\begin{enumerate}
 \item the limits 
 $$\varphi^{\pm}(x)=\lim_{n\to\infty}\frac{1}{n}\sum_{k=0}^{n-1}\varphi(f^{\pm k}(x))$$
exist $m$-almost every $x\in M$.
\item We have $$\varphi^{+}(x)=\varphi^{-}(x)\qquad m\text{-almost every }x$$
\item The functions $\varphi^{\pm}$ are in $L^{1}_{m}(M)$ (that is, they are integrable), and 
$$\int \varphi^{\pm}dm=\int \varphi dm.$$
\end{enumerate} 
\end{theorem}
This theorem provides an alternative way to establish ergodicity:
\begin{corollary}\label{birkhoff.ergodicity}
 Let $f:M\to M$ be a diffeomorphism preserving the measure $m$. Then $f$ is ergodic if and only if 
 for all $\varphi\in C^{0}(M)$, 
\begin{equation}\label{eq.birkhoff.ergodic}
 \varphi^{\pm}(x)=\int\varphi dm\qquad m\text{-almost every }x\in M
 \end{equation}
\end{corollary}
\begin{remark}\label{remark.birkhoff}
 For future use, to establish ergodicity it is enough to check Equation \eqref{eq.birkhoff.ergodic} just for an $L^{1}_{m}$-dense family of continuous functions $\cD$. Actually it is enough to check that each $\varphi\in\cD$ satisfies that $\varphi^{+}$ (or $\varphi^{-}$) is constant $m$-almost everywhere.  
\end{remark}
\subsection{Exercises}\label{constant.manifold}
\begin{enumerate}
 \item Prove that $\varphi^{+}$ is constant on $W^{-}$ (also on $W^{s}$) for all continuous maps $\varphi$. Analogously, $\varphi^{-}$ is constant on $W^{+}$. 
\end{enumerate}
\section{The set of typical points}
In this section, let $f:M\to M$ be a $C^{1+\alpha}$ diffeomorphism preserving the volume measure $m$. 
Then we have the following:
\begin{proposition}[Typical points] \label{typical} There exists an invariant set $\cS$, the set of {\em typical points} with $m(\cS)=1$ such that for all $\varphi\in \cD$, where $\cD\subset C^{0}(M)$ is an $L^{1}_{m}(M)$ dense set,  all $x_{0}\in \cS$ satisfy
\begin{enumerate}
 \item $\varphi^{+}(x_{0})=\varphi^{+}(w)$ for all $w\in W^{-}(x_{0})$, 
 \item $\varphi^{+}(x_{0})=\varphi^{+}(w)$ for $m^{+}_{x_{0}}$-almost every $w$ in $W^{+}(x_{0})$. 
\end{enumerate}
Here $m^{+}_{x_{0}}$ denotes the conditional measures induced by a measurable partition subordinated to the unstable partition $W^{+}$. 
\end{proposition}
The existence of the set of typical points $\cS$ follows essentially from the Birkhoff ergodic theorem in the previous section, together with the absolute continuity properties discussed in last chapter. Let us give a sketch of the proof. \par
Item (1) follows from Exercise (1) in Subsection \ref{constant.manifold}.\par
Now, from the Birkhoff Ergodic Theorem it follows the existence of a set $A$ with $m(A)=1$ such that 
$$\varphi^{+}(x)=\varphi^{-} (x)\qquad\forall x\in A.$$
Since both $\varphi^{+}$ and $\varphi^{-}$ are invariant functions, that is $\varphi^{\pm}\circ f(x)=\varphi^{\pm}(x)$ for all $x\in A$, then $A$ is an invariant set. Since $m^{+}_{x}$ is the system of conditional measures subordinated to $W^{+}$, the Fubini-type Equation \eqref{fubini} yields
$$m^{x}(W^{+}(x)\setminus A)=0\qquad m\text{-almost every }x\in M$$
Then, for $m$-almost all $x\in A$, $m^{+}_{x}$-almost every $w$ in $W^{+}(x)$, we will have:
$$\varphi^{+}(x)=\varphi^{-}(x)=\varphi^{-}(w)=\varphi^{+}(w).$$
Pay attention to the fact that the invariant set $A$ satisfies both claims but only for the function $\varphi$. Now take a countable family $\cD$ of continuous functions that are $L^{1}_{m}(M)$-dense in $L^{1}_{m}$. The interseccion of all $A_{\varphi}$ with $\varphi\in\cD$ provides us with the desired set $\cS$ of typical points. \ep
\section{The Hopf argument}\label{hopf}
We will use the Hopf argument to prove a classical theorem:
\begin{theorem}[Anosov-Sinai 1967] Let $f:M\to M$ be a $C^{1+\alpha}$ Anosov diffeomorphism (that is, a diffeomorphism such that the whole manifold $M$ is a hyperbolic set for $f$), preserving a volume measure $m$. Then $f$ is ergodic.
\end{theorem}
Let us take a countable family $\cD\subset C^{0}(M)$ dense in $L^{1}_{m}(M)$ and the corresponding set of typical points $\cS$, as obtained in Proposition \ref{typical}. Let $\varphi\in \cD$. Our strategy will be to see that $\varphi^{+}$ is locally constant modulo a Lebesgue zero set. From compactness of $M$ and Remark \ref{remark.birkhoff} the claim follows. So take $x_{0},y_{0}\in \cS$ that are really close, say $d(x_{0},y_{0})<\eps$. If $\eps>0$ is sufficiently small, for all $x\in B_{\eps}(x_{0})$, the local stable manifold of $x$ will cut both the local unstable manifold of $x_{0}$ and of $y_{0}$. This follows from the Stable Manifold Theorem and the transversality of $E^{s}$ and $E^{u}$. \par
Since $x_{0}\in \cS$ is a typical point, $m^{+}_{x_{0}}$-almost every $w\in W^{+}_{loc}(x_{0})$ will satisfy $\varphi^{+}(w)=\varphi^{+}(x_{0})$. See Figure. 
\begin{figure}[h]
\begin{center}
 \includegraphics[width=\textwidth]{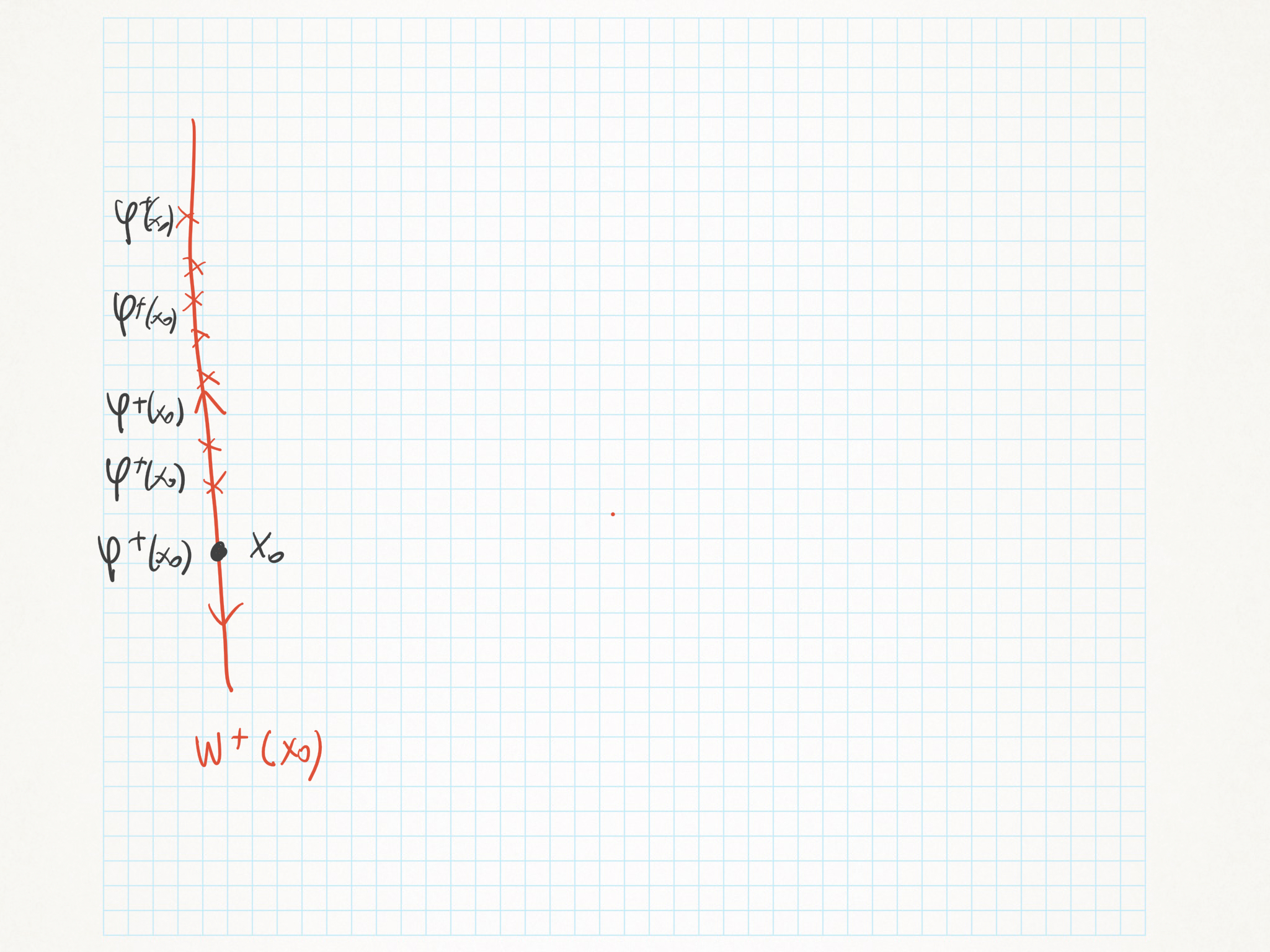}
\end{center}
\end{figure}
Call $A=\{w\in W^{+}_{loc}(x_{0}):\varphi^{+}(w)=\varphi^{+}(x_{0})\}$. Now consider the stable holonomy map $h$ between $W^{+}_{loc}(x_{0})$ and $W^{+}_{loc}(y_{0})$. Since there is transverse absolute continuity, $m^{+}_{y_{0}}(h(A))>0$. But $\varphi^{+}$ is constant on $W^{-}$-leaves, so for all points $w\in h(A)\subset W^{+}_{loc}(y_{0})$, $\varphi^{+}(w)=\varphi^{+}(x_{0})$. See Figure. 
\begin{figure}[h]
\begin{center}
 \includegraphics[width=\textwidth]{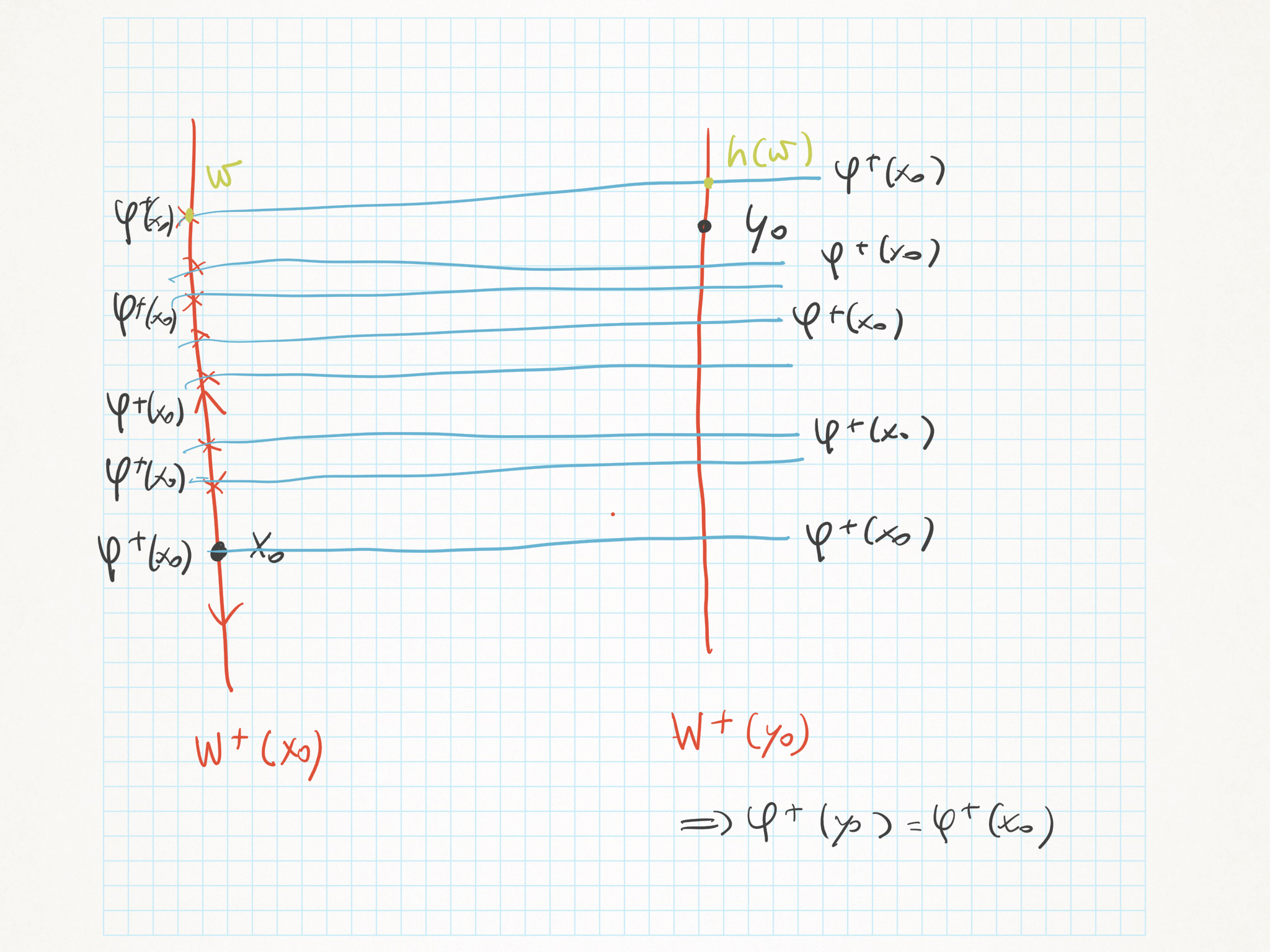}
\end{center}
\end{figure}
Now we have an $m^{+}_{y_{0}}$-positive measure set in $W^{+}_{loc}(y_{0})$ such that $\varphi^{+}$ takes the value $\varphi^{+}(x_{0})$ on it. Since $y_{0}\in \cS$ is a typical point, $\varphi^{+}(y_{0})$ has to be equal to $\varphi^{+}(x_{0})$. And so we have proved that $\varphi^{+}$ is almost everywhere locally constant. Since $\cD$ is a dense set in $L^{1}_{m}(M)$, we have proved that $f$ is ergodic. \ep
\chapter{An ergodic spectral decomposition theorem}
\section{Smale's Spectral Decomposition Theorem}
In the hyperbolic world, we have an interesting tool to codify the dynamics when periodic points are abundant. Concretely, if $f:M\to M$ is a diffeomorphism, we may define the {\em non-wandering set of $f$} by
$$NW(f)=\{x\in M:\forall U\ni x\text{ open}, \exists n\ne0 \text{ s.t. }f^{n}(U)\cap U\ne\emptyset\}$$
We will say that $f$ is {\em Axiom A} if 
\begin{enumerate}
 \item $NW(f)$ is a hyperbolic set, 
 \item $NW(f)=\overline{\per(f)}$. 
\end{enumerate}
$\per(f)$ denotes the set of periodic points, and $\per_{H}(f)$ will denote the set of hyperbolic periodic points. \par
Smale has given a pretty good description of the dynamics in the non-wandering set of an Axiom A diffeomorphism. For a nice overview of this theorem and uniform hyperbolicity in general, I recommend the books \cite{hasselblattkatok} and \cite{brinstuck}. 
\begin{theorem}[Smale's Spectral decomposition theorem] \label{smale}Let $f:M\to M$ be an Axiom A diffeomorphism. Then there exist compact invariant disjoint sets $\Lambda_{1},\dots,\Lambda_{n}$, called {\em basic sets} satisfying
$$NW(f)=\Lambda_{1}\cup\dots\cup\Lambda_{n}.$$
 $f$ is transitive on each $\Lambda_{i}$. Each of these sets $\Lambda_{i}$ decomposes as a finite union of compact disjoint sets:
 $$\Lambda_{i}=\Gamma_{i}^{1}\cup\dots\cup\Gamma^{j_{i}}_{i}, $$
 so that $f(\Gamma^{k}_{i})=\Gamma^{k+1}_{i}$ for all $k\in \Z\mod j_{i}$. We also have
 $$f^{j_{i}}|_{\Gamma^{k}_{i}}\quad\mbox{
\begin{tabular}{l}
\text{ is conjugated to a shift of finite type, in particular, }\\ \text{ it is topologically mixing}.
\end{tabular}}
$$
\end{theorem}
Remember that $f$ is {\em transitive} on an invariant set $\Lambda$ if there is a point $x$ in $\Lambda$ with a dense orbit: $\overline{o(x)}=\Lambda$. $f$ is {\em topologically mixing} on an invariant set $\Gamma$ if for any two non-empty open sets $U,V$ there exists $N>0$ such that 
$$f^{n}(U)\cap V\ne\emptyset \quad \forall n\geq N.$$
(Sorry, we are not defining shifts of finite type for the time being, but in a future version -if there is one- there could be in an Apendix about them or a separate chapter. Meanwhile, you can check the excellent books of Hasselblatt-Katok \cite{hasselblattkatok} and Brin-Stuck \cite{brinstuck}, for instance)
\subsection{Exercises}
\begin{enumerate}
 \item Check that $NW(f)$ is always compact, invariant and non-empty. 
 \item Check that topologically mixing implies transitivity.
\end{enumerate}
\section{Homoclinic classes}\label{homoclinic}
Even though we are not providing any proof of Smale's Spectral Decomposition Theorem \ref{smale}, we can improve the description of the sets $\Lambda_{i}$ and $\Gamma^{k}_{i}$. 
Given two hyperbolic periodic points $p,q\in\per_{H}(f)$, we will denote $p\approx q$ if
$$W^{s}(q)\transv W^{u}(p)\ne\emptyset\quad\text{ and }\quad W^{u}(q)\transv W^{s}(p)\ne\emptyset$$
\begin{figure}[h]
\begin{center}
 \includegraphics[width=\textwidth]{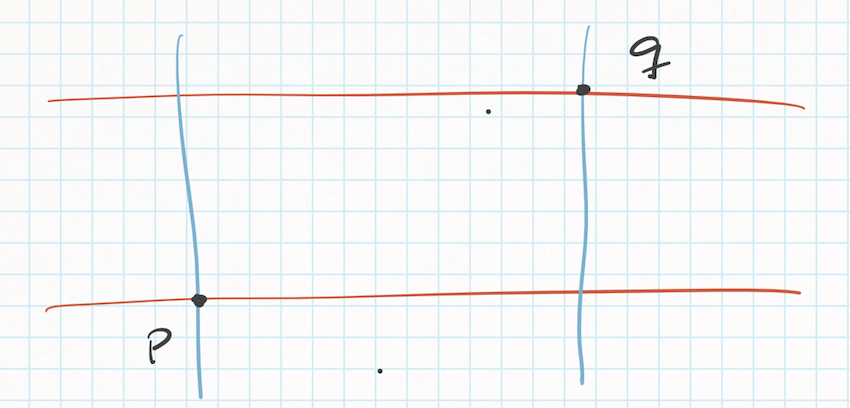}
\end{center}
\end{figure}
We will say that two hyperbolic periodic points $p,q\in\per_{H}(f)$ are {\em homoclinically related} if $q\approx f^{k}(q)$ for some $k\in \Z$, and we denote it $q\in H(p)$. The {\em homoclinic class} of $p$ is defined by $\overline{H(p)}$. Then all basic sets $\Lambda_{i}$ are homoclinic classes. Namely, for each $i=1,\dots,n$ there exists $p_{i}\in\per_{H}(f)$ such that 
$$\Lambda_{i}=\overline{H(p_{i})}. $$
We can also describe the sets $\Gamma^{k}_{i}$:
$$\Gamma_{i}^{k}=\overline{\{q:q\approx f^{k}(p_{i})\}}.$$
\section{Pesin's ergodic spectral decomposition theorem}
In 1977, Pesin provided a spectral decomposition theorem in the {\em non-uniformly hyperbolic set}, a.k.a. {\em Pesin region}:
$$\nuh(f)=\{x\in\cR:\dim E^{0}_{x}=0\}$$
\begin{theorem}[Pesin's Ergodic Spectral Decomposition Theorem \cite{pesin77}] Let $f:M\to M$ be a $C^{1+\alpha}$-diffeomorphism preserving a volume probability measure $m$. Suppose $m(\nuh(f))>0$. Then
there exists a countable family of $\mod 0$-pairwise disjoint invariant sets $\{\Lambda_{n}\}_{n\in\N}$, that is, $m(\Lambda_{i}\cap\Lambda_{j})=0$ if $i\ne j$ such that
$$\nuh(f)=\Lambda_{1}\cup\dots\cup \Lambda_{i}\cup\dots$$
where $f|_{\Lambda_{i}}$ is ergodic for all $i\in\N$. $\Lambda_{i}$ are the {\em ergodic components}.\par Moreover, each $\Lambda_{i}$ decomposes as a finite union of $\mod 0$-disjoint sets:
$$\Lambda_{i}=\Gamma_{i}^{1}\cup\dots\cup \Gamma^{j_{i}}_{i}$$
 such that $f(\Gamma^{k}_{i})=\Gamma^{k+1}_{i}$ for all $k\in\Z\mod j_{i}$. And
 $$f^{j_{i}}|_{\Gamma^{k}_{i}}\quad\mbox{
\begin{tabular}{l}
 isomorphic to a Bernoulli shift, in particular, mixing. 
\end{tabular}
}$$
\end{theorem}
We would like to complete the description of the ergodic components as an analog of the homoclinic classes described in Section \ref{homoclinic}, and that is what we are going to do next. 
\section{Ergodic homoclinic classes}
In 2011, Hertz, Hertz, Tahzibi and Ures introduced the ergodic homoclinic classes \cite{HHTU11}. 
These sets are defined as follows: given $p\in\per_{H}(f)$ the {\em stable ergodic homoclinic class of $p$} is the set 
$$\ehc^{-}(p)=\{x\in\cR: W^{-}(x)\transv W^{u}(o(p))\ne\emptyset\}$$
$o(p)$ denotes the orbit of $p$. This set is invariant and {\em saturated by Pesin stable manifolds}. This means that if $x\in\ehc^{-}(p)$ then the whole Pesin stable manifold $W^{-}(x)\subset \ehc^{-}(p)$. This follows directly from the definition. 
Analogously, we define the {\em unstable ergodic homoclinic class of $p$} by
$$\ehc^{+}(p)=\{x\in\cR: W^{+}(x)\transv W^{s}(o(p))\ne\emptyset\}$$
$\ehc^{+}(p)$ is invariant and saturated by Pesin unstable manifolds; that is, if $x\in \ehc^{+}(p)$, then the whole Pesin unstable manifold $W^{+}(x)\subset \ehc^{+}(p)$. 
\begin{figure}
\centering
\begin{subfigure}{.5\textwidth}
  \centering
  \includegraphics[width=1\linewidth]{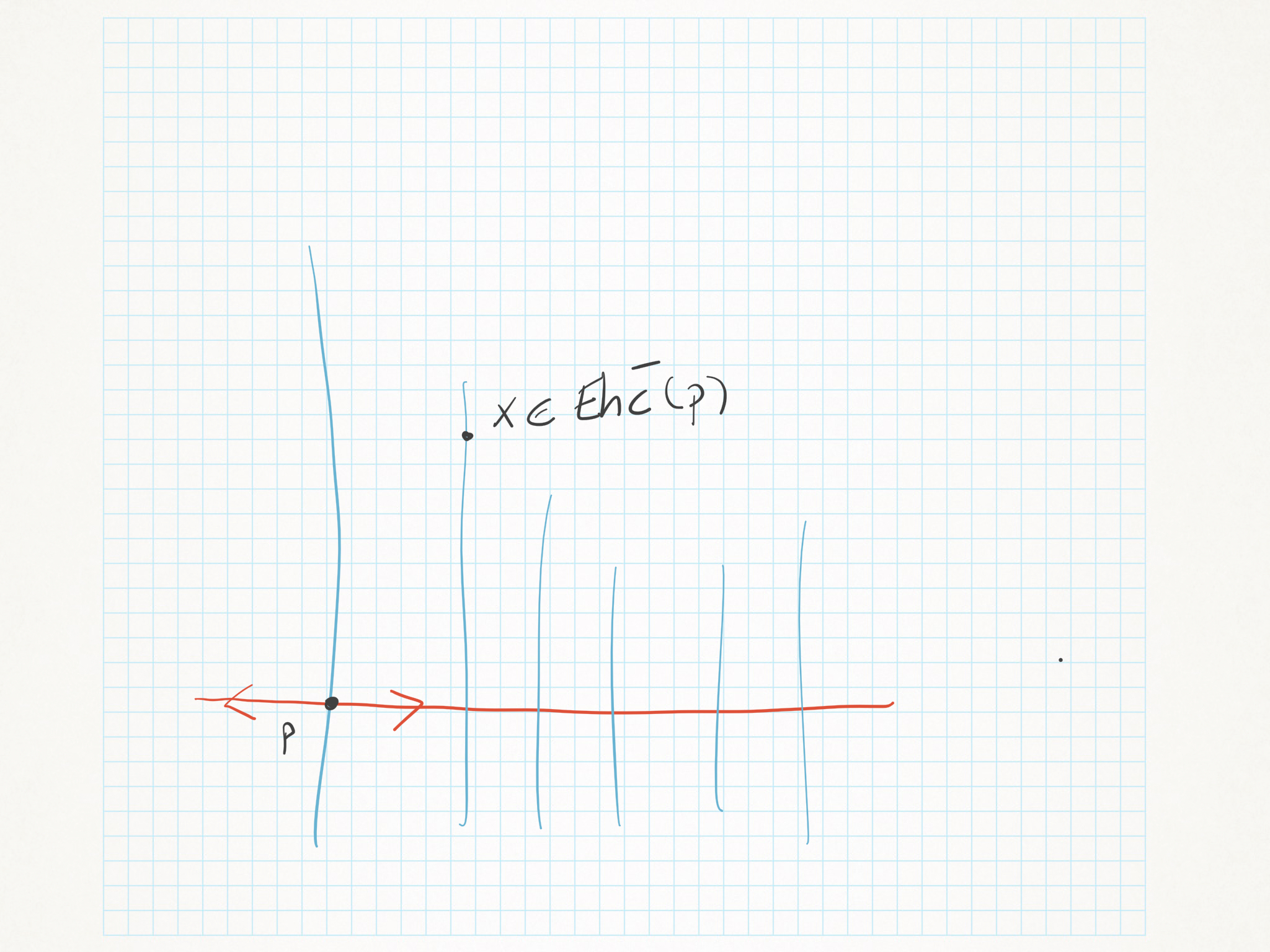}
  \caption{$\ehc^{-}(p)$}
\end{subfigure}%
\begin{subfigure}{.5\textwidth}
  \centering
  \includegraphics[width=1\linewidth]{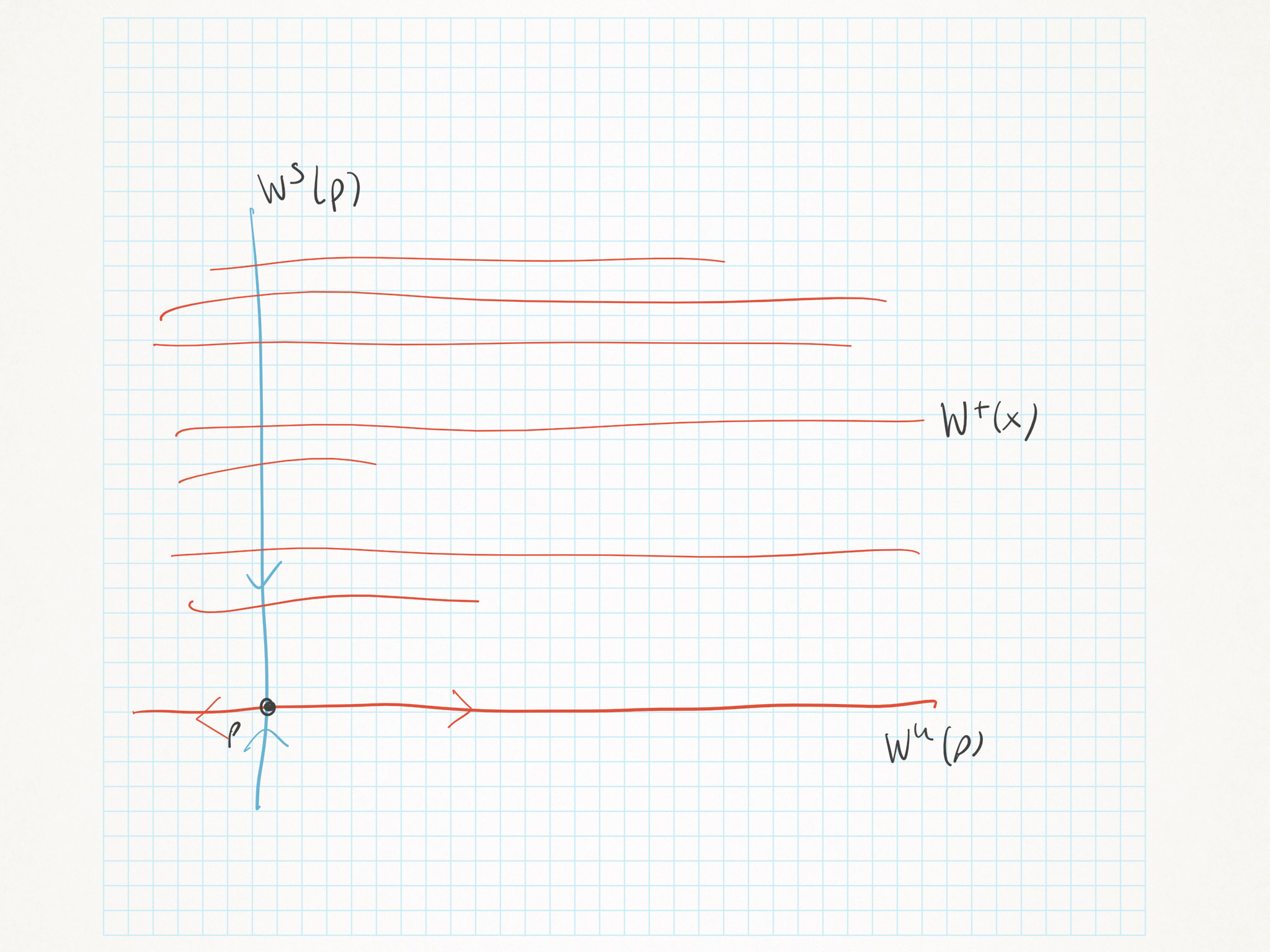}
  \caption{$\ehc^{+}(p)$}
\end{subfigure}
\caption{Ergodic homoclinic classes}
\end{figure}

We define the {\em ergodic homoclinic class of $p$} by
$$\ehc(p)=\ehc^{-}(p)\cap \ehc^{+}(p)$$
We have the following:
\begin{theorem}[Criterion HHTU 2011] Let $f:M\to M$ be a $C^{1+\alpha}$ diffeomorphism preserving a volume measure $m$. Let $p\in\per_{H}(f)$. If
$$m(\ehc^{-}(p))>0\quad\text{ and }\quad m(\ehc^{+}(p))>0, $$
 then
\begin{enumerate}
 \item \label{item1}$\ehc^{-}(p)\zeroeq \ehc^{+}(p)\zeroeq\ehc(p),$
 \item \label{item2}$f|_{\ehc(p)}$ is ergodic
 \item \label{item3}$\ehc(p)\subset\nuh(f).$
\end{enumerate}
\end{theorem}
In the theorem above $\zeroeq$ means the sets coincide modulo a zero Lebesgue measure set. \par
Note that Item (\ref{item3}) follows directly from the definition of $\ehc(p)$. Indeed, let $s=\dim W^{s}(p)$ and $u=\dim W^{u}(p)$. Since $p$ is hyperbolic $s+u=n=\dim M$. 
If $x\in \ehc(p)$ then $x\in\ehc^{-}(p)$, so by definition $W^{-}(x)\transv W^{u}(o(p))\ne\emptyset$. This implies that $\dim E^{-}_{x}=\dim W^{-}(x)\geq n-u =s$. Analogously, $x\in\ehc^{+}(p)$ so $W^{+}(x)\transv W^{s}(o(p))\ne\emptyset$. This implies $\dim E^{+}_{x}=W^{+}(x)\geq n-s=u$. Since $E^{+}_{x}\cap E^{-}_{x}=\{0\}$, we have $\dim E^{0}_{x}=\{0\}$ and therefore $x\in \nuh(f)$. \par
The proof of Item (\ref{item2}) follows the Hopf argument, which we described in Section \ref{hopf}. We will fill the details here for completeness. Let us begin by recalling a folkloric lemma from the hyperbolic world, the {\em inclination lemma} a.k.a. the {\em $\lambda$-lemma}. Let $p\in\per_{H}(f)$ (for simplicity we will assume it to be fixed). Then for any disc $D$ transversely intersecting $W^{s}(p)$, and any disc $B\subset W^{u}(p)$ for all sufficiently large $n>0$ there exists a disc $D_{B}\subset f^{n}(D)$ that is $C^{1}$-close to $B$. We will not formalize this concept more than this, but see the picture to get an idea. 
\begin{figure}[h]
\begin{center}
 \includegraphics[width=\textwidth]{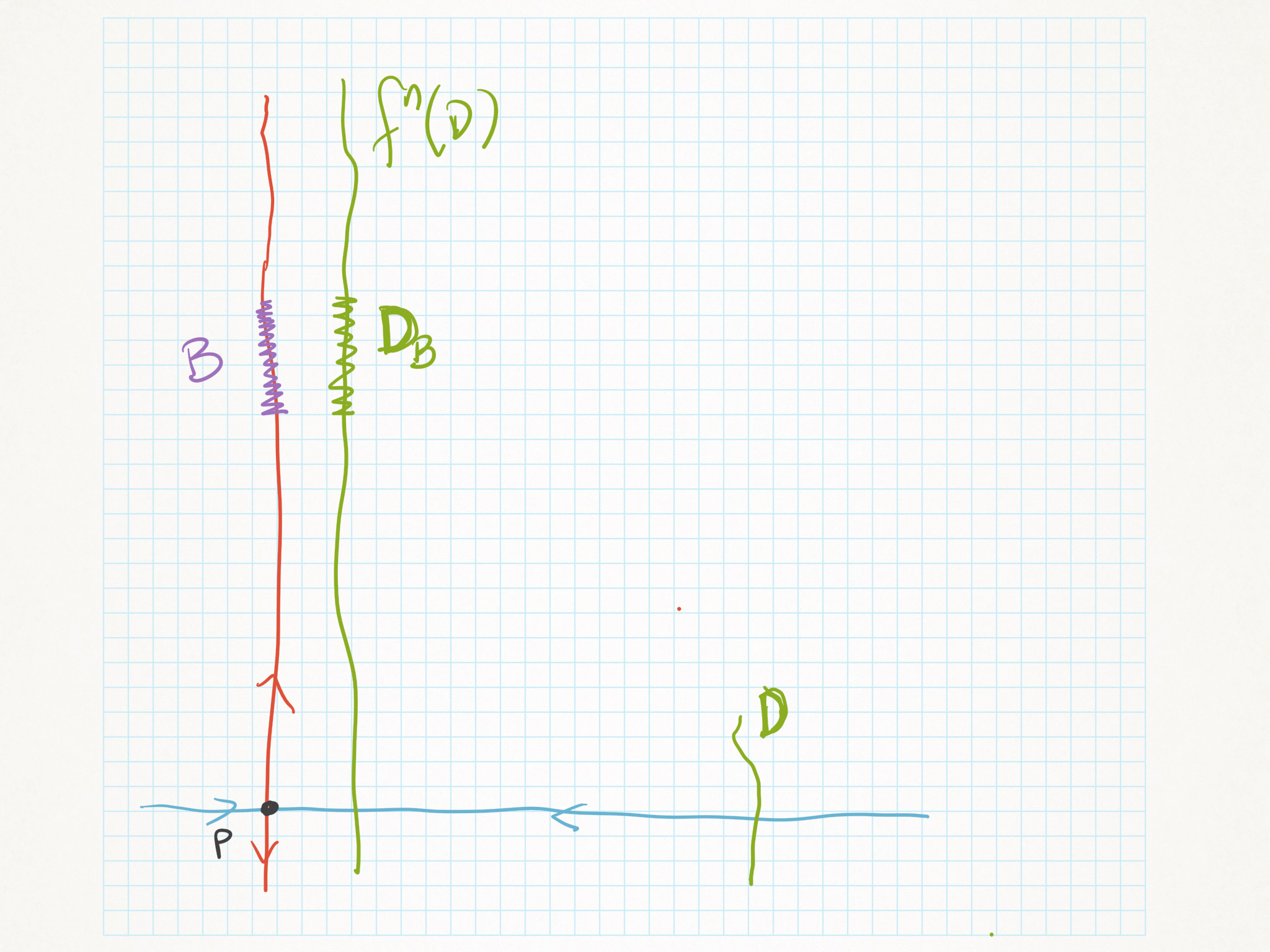}
 \caption{The inclination lemma or $\lambda$-lemma}
\end{center} 
\end{figure}\par
Let $\cS$ be the set of typical points obtained in Proposition \ref{typical}, and let $\varphi\in\cD$ a continuous function, where $\cD$ is a countable $L^{1}_{m}$-dense set as in Proposition \ref{typical}. Let $x_{0},y_{0}$ be Lebesgue density points of  $\cS\cap\ehc(p)$. We want to see that $\varphi_{+}(x_{0})=\varphi_{+}(y_{0})$. From this it will follow that $f$ is ergodic when restricted to $\ehc(p)$. Since $y_{0}\in\cS$ is a typical point, we have that $\varphi_{+}(w)=\varphi^{+}(y_{0})$ for $m^{+}_{y_{0}}$-almost every $w\in W^{+}(y_{0})$. Since $x_{0}$ is a Lebesgue density point of $\cS\cap\ehc^{-}(p)$, and this is an invariant set, we may assume $x_{0}$ is very close to $W^{u}(p)$. Indeed, $W^{-}(x_{0})\transv W^{u}(o(p))\ne\emptyset$, so by taking an arbitrarily large iterate of $f$, we will get $f^{n}(x_{0})$ super close to $W^{u}(p)$. This does not affect the value of $\varphi^{+}$, since $\varphi^{+}(x_{0})=\varphi^{+}(f^{n}(x_{0}))$. Let us say $x_{0}$ is $\delta$-close to $W^{u}(p)$, with very small $\delta>0$. Since $m^{+}_{x}$ is a canonical set of conditional measures subordinate to $W^{+}$, for almost every $x\in B_{\delta}(x_{0})$, $m^{+}_{x}(\ehc(p)\cap W^{+}(x))>0$. \par
Now, since $y_{0}\in\ehc(p)\subset \ehc^{+}(p)$, $W^{+}(y_{0})\transv W^{s}(o(p))\ne\emptyset$. The inclination lemma then implies that $f^{n}(W^{+}(y_{0}))$ will contain a disc $D$ inside $B_{\delta}(x_{0})$ which is $C^{1}$-close to a disc in $W^{u}(p)$ for sufficiently large $n>0$. The transverse absolute continuity property then implies that the stable holonomy map takes the set of points $w$ in $D$ such that $\varphi^{+}(w)=\varphi^{+}(y_{0})$, which has $m^{+}_{y}$- positive measure into a set of $m^{+}_{x_{0}}$- positive measure of $W^{+}(x_{0})$. Since $\varphi^{+}$ is constant on $W^{-}$-leaves, this implies that $W^{+}(x_{0})$ contains a set of positive $m^{+}_{x_{0}}$-measure such that all $w$ in this set satisfy $\varphi^{+}(w)=\varphi^{+}(y_{0})$. Since $x_{0}\in\cS$ is a typical point, this implies that $\varphi^{+}(x_{0})=\varphi^{+}(y_{0})$, and we have proved Item (\ref{item2}). The figure can perhaps help understand the argument. 
\begin{figure}[h]
\begin{center}
 \includegraphics[width=\textwidth]{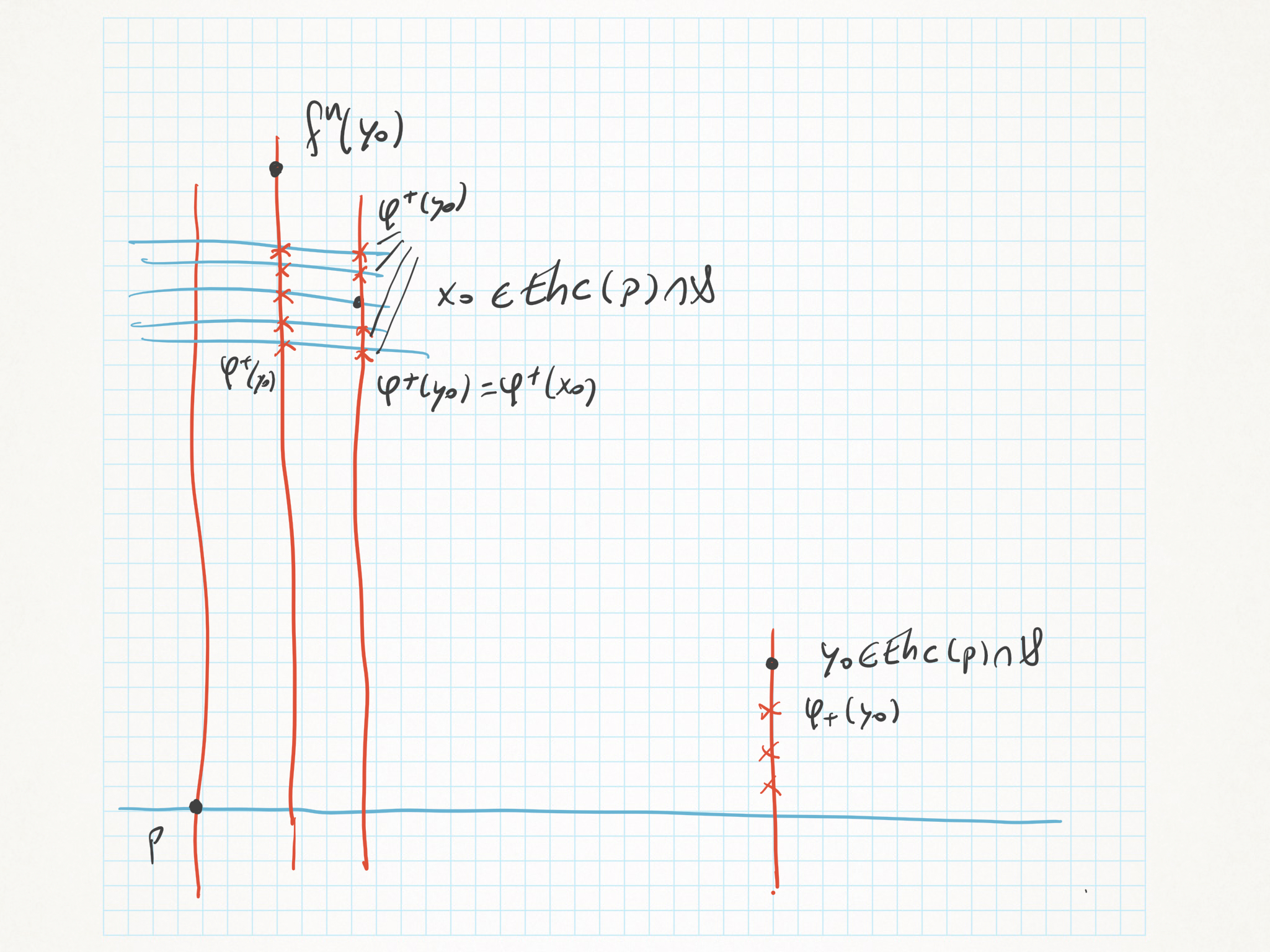}
 \caption{A proof that $f$ is ergodic on $\ehc(p)$}
\end{center} 
\end{figure}\par
Now we only have to see Item (\ref{item1}), which is the most involved one. From Proposition \ref{typical} it can be deduced the existence of an invariant set $\cT$ with $m(\cT)=1$, such that for all $x\in\cT$
\begin{equation}\label{typical2}
 \mathbbm{1}_{\pm}(x)=\mathbbm{1}_{\pm}(w)\quad m^{\pm}_{x}\text{-almost every } w\in W^{\pm}(x),
\end{equation}
 where $\mathbbm{1}_{\pm}$ is the characteristic function of the set $\ehc^{\pm}(p)$, that is, the function that is equal to one on $\ehc^{\pm}(p)$ and equal to zero otherwise. \par
Why this is true is an act of faith, but if you are curious, you can take a look at Lemma 4.3 in HHTU paper \cite{HHTU11}. \par
Take $x_{0}$ a Lebesgue density point of $\ehc^{-}(p)\cap\cT$. We want to prove that $x_{0}\in\ehc^{+}(p)$. Now take $y_{0}$ a Lebesgue density point of $\ehc^{+}(p)\cap\cT$. See figure. 
\begin{figure}[h]
\begin{center}
 \includegraphics[width=\textwidth]{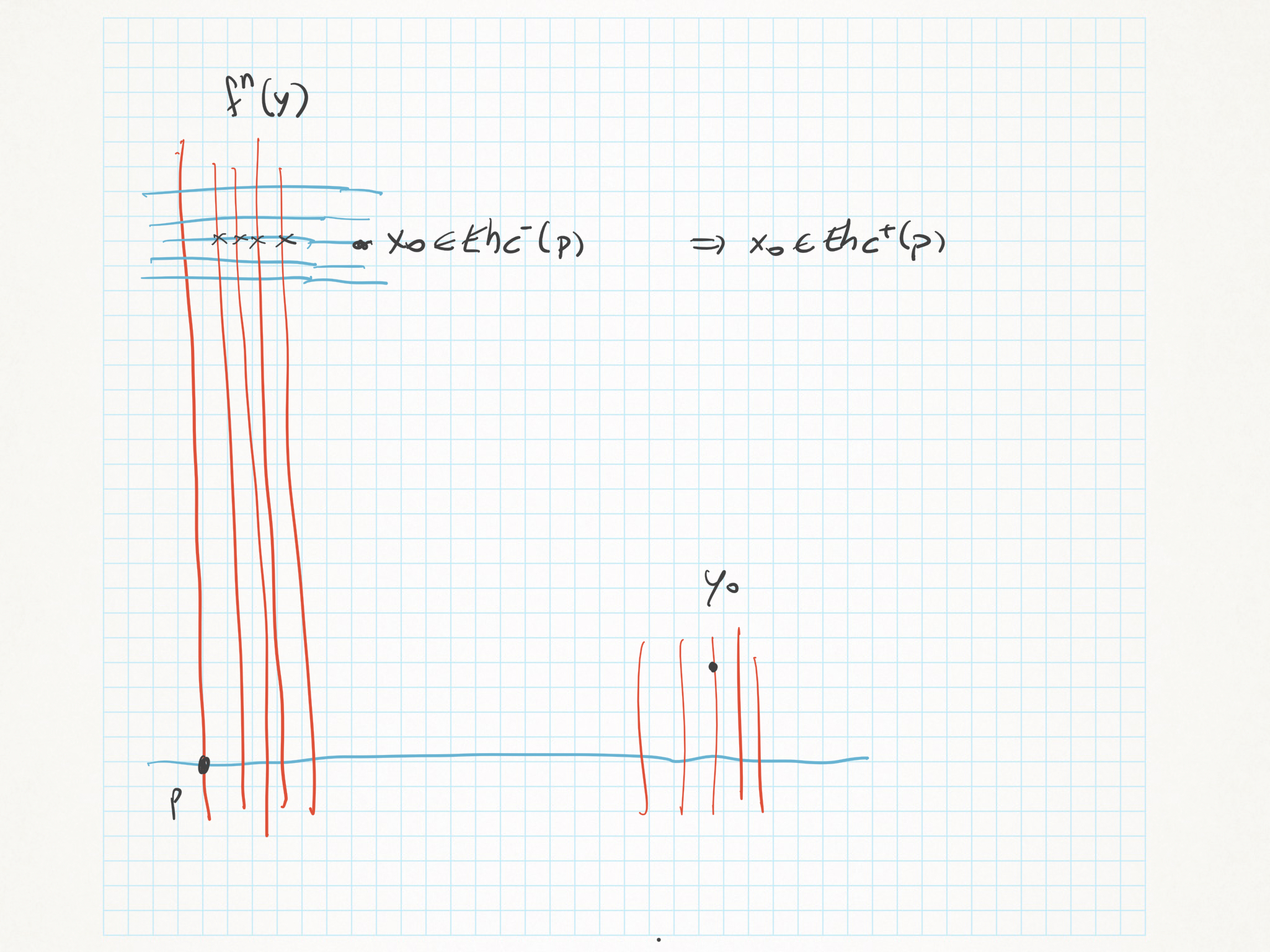}
 \caption{A proof of the criterion HHTU}
\end{center} 
\end{figure}\par
Then the inclination lemma implies that for all $y\in B_{\eps}(y_{0})\cap\ehc^{+}(p)\cap\cT$ there is $n_{y}>0$ such that $f^{n_{y}}(W^{+}(y))\transv W^{-}(x_{0})\ne\emptyset$. Transverse absolute continuity and the countability of $n$'s imply there exists an $n>0$ for which the set of $w$ such that 
$f^{n}(W^{+}(w))\transv W^{-}(x_{0})\ne\emptyset$ has positive measure in $W^{-}(x_{0})$. This implies that
$$m^{-}_{x_{0}}(w:\mathbbm{1}_{+}(w)=1)>0,$$
then our choice of $\cT$ (see Equation \eqref{typical2}) implies $x_{0}\in\ehc^{+}(p)$, and with this we have finished the proof of the criterion. \ep
\begin{remark} Hey, not so fast! Ok, I know you have liked the proof, but I've cheated a bit in this last piece. What happens if $\dim W^{-}(x_{0})>s=\dim W^{s}(p)$? Well, in that case you can easily sub-foliate $W^{-}(x_{0})$ by a measurable sub-foliation of dimension $s$, let us call it $W^{--}$. This sub-foliation has its canonical system of conditional measures $m^{--}_{x}$ with $x\in W^{-}(x_{0})$. The previous argument implies that 
$$m^{--}_{x}(w:\mathbbm{1}_{+}(w)=1)>0\qquad \text{for a }m^{-}_{x_{0}}-\text{positive measure set} $$
 The fact that $W^{--}$ is a measurable sub-foliation then implies that 
 $$m^{-}_{x_{0}}(w:\mathbbm{1}_{+}(w)=1)>0,$$
 and the proof follows as before. A fortiori this situation will happen only in a zero Lebesgue measure set. 
\end{remark}
\section{Pesin's ergodic spectral decomposition theorem revisited}
The criterion HHTU described in previous section gives us a geometrical description of ergodic components. But hey! we need a hyperbolic periodic point. This was done by Katok. 
\begin{theorem}[Katok's closing lemma \cite{katok80}] Let $f:M\to M$ be a $C^{1+\alpha}$-diffeomorphism and $\mu$ a hyperbolic probability measure. Then:
\begin{enumerate}
\item $\supp(\mu)\subset\overline{\per_{H}(f)}$, and
\item $\mu$-almost every $x\in M$ there exists $p\in \per_{H}(f)$ such that 
$$x\in\ehc(p).$$ 
\end{enumerate}
\end{theorem}
Of course Katok did not formulate it like that, but he indeed proved it, and it is all we need now. 
With Katok's closing lemma and the HHTU criterion, we are now able to reformulate Pesin's ergodic spectral decomposition theorem:
\begin{theorem} [Pesin's ergodic spectral decomposition theorem revisited \cite{pesin77}, \cite{katok80}, \cite{HHTU11}] Let $f:M\to M$ be a $C^{1+\alpha}$-diffeomorphism preserving a volume probability measure $m$. Suppose $m(\nuh(f))>0$. Then
there exists a countable family of hyperbolic periodic points $\{p_{n}\}_{n\in\N}\subset\per_{H}(f)$ such that 
$$\nuh(f)=\ehc(p_{1})\cup\dots\cup\ehc(p_{n})\cup\dots$$
where $\ehc(p_{n})$ is the ergodic homoclinic class of $p_{n}$, the ergodic homoclinic classes are two-by-two disjoint modulo a zero Lebesgue measure set. 
\par Each ergodic homoclinic class can be written as 
$$\ehc(p_{i})=\Gamma_{i}^{1}\cup\dots\cup \Gamma^{j_{i}}_{i}$$
where $\Gamma{i}^{k}=\{x: W^{-}(x)\transv W^{u}(f^{k}(p))\quad\text{and}\quad W^{+}(x)\transv W^{s}(f^{k}(p))\}$, and
 $$f^{j_{i}}|_{\Gamma^{k}_{i}}\quad\mbox{
\begin{tabular}{l}
 isomorphic to a Bernoulli shift, in particular, mixing. 
\end{tabular}
}$$
\end{theorem}

\end{document}